# The *khôra*
# and the two-triangle universe of Plato's *Timaeus*

Luc Brisson
Centre Jean Pépin
CNRS-ENS/UMR 8230

Salomon Ofman
Institut mathématique de Jussieu-
Paris Rive Gauche/HSM
Sorbonne Université-Université
Paris Sorbonne cité

***Abstract.*** The main purpose of this article is to try to understand the connection between the physical universe and the mathematical principles that underlies the cosmological account of the *Timaeus*. Aristotle's common criticism of Plato's cosmology is that he confuses mathematical and physical constructions. Namely, the *Timaeus* is the first cosmology founded on mathematical physics. We give a new translation of *Timaeus* 31b-32b, an important passage to understand the connection between mathematics and physics in Timaeus' cosmological construction.

This article is the first of a series about the *khôra*. We will restrict our focus here to the much-debated question of the primary elements in the *khôra*, the components of the whole physical world reduced, in an extraordinary elegant construction, to two right triangles.

## 1. Introduction

The *Timaeus* is Plato's only cosmological dialogue where one of the most difficult questions of his doctrine is considered in some detail: the so-called 'participation' of 'sensible world' in the intelligible world, or the Forms or Ideas. The main interlocutor is Timaeus, who is presented as coming from Locri, a city in the south of Italy, where the Pythagoreans were active,[1] someone who 'knows more of astronomy than the rest of us and has made knowledge of the nature of the universe his chief object'.[2] Then begins a long cosmological account based on a new 'genus' ('*genos*'), the '*khôra*',[3] which has not appeared and will never appear again in any other dialogue. Its usual meaning in Greek is 'space'; however, Timaeus emphasizes that it is something more. The *khôra* is not only the place where perpetual changes of the sensible bodies occur, but also the unqualified (49c7-50a4) and unchanging (52a8-d1) 'genus' required to explain these perpetual changes; it is the dynamical whole,[4] consisting of these bodies in perpetual change, their 'nurse' ('*tithênên*'), their nourisher

---

[1] Thus, it is usually assumed that Plato wanted to present him as a Pythagorean (as already, Proclus in his *Commentary on the* Timaeus (Runia-Share (2008), II, 223.5-7)). However, one has to be careful since there is no other source of information about Timaeus than Plato's eponymous dialogue (cf. Nails (2002), p. 293); on the subsequent imaginative development of the personage, see C. Macris, *DPhA* VII 2018, 987-1009.

[2] 27a3-4.

[3] As Zeyl remarks, the only 'unambiguous denomination of the receptacle' found in the *Timaeus* is '*khôra*' (Zeyl (2010), p. 118). Nevertheless, this 'third genus (…) which exists always and cannot be destroyed' (52a8), is given several different metonymic names throughout the dialogue, the 'receptacle' ('*hypodokhên*', 49a6) being one among others (cf. *infra*, note 5).

[4] The dynamical aspect of the *khôra* will be studied in a next article.

('*trophon*') and their 'mother' ('*mêtera*').[5] This difficult, even puzzling, 'new genus'[6] has been the subject of many different interpretations almost from Plato's time until our days.[7] However, pro or con, these were usually written through an Aristotelian lens,[8] and it was soon assimilated to Aristotelian matter (the '*hylê*'). We have chosen another approach here, avoiding substituting Plato's concepts for Aristotelian ones. From this perspective, Plato seems paradoxically to agree, in anticipation, with most of Aristotle's criticisms, except that he does not consider them as such.[9]

## 2. The binding proportions

Timaeus' narrative about the 'nature of the universe'[10] begins when he agrees to give an account of it from its origins to the emergence of men.[11] His first claim concerns the theory of the four elements: earth, water, air, fire, as components of all sensible bodies.[12] Such a claim is already found in Empedocles' texts under the name of the four roots (*rhizai*), largely predating Timaeus', i.e. Socrates' time.[13] However, Timaeus immediately proposes a complete novelty with respect to his predecessors:[14] the four elements are not the ultimate components, and the theory is not based on some physical or experimental reasoning, but on a mathematical one.

### i) *Timaeus* (31b-32b) and its translation

| | |
|---|---|
| Now that which comes to be, must be bodily form, and so visible and tangible; and nothing could ever become visible without fire, nor tangible without something solid, nor solid without earth. That is why, when he began to put the body of the world together [*synistanai*], the god came to make it out of fire and earth.<br>But it isn't possible to put together | σωματοειδὲς δὲ δὴ καὶ ὁρατὸν ἁπτόν τε δεῖ τὸ γενόμενον εἶναι, χωρισθὲν δὲ πυρὸς οὐδὲν ἄν ποτε ὁρατὸν γένοιτο, οὐδὲ ἁπτὸν ἄνευ τινὸς στερεοῦ, στερεὸν δὲ οὐκ ἄνευ γῆς: ὅθεν ἐκ πυρὸς καὶ γῆς τὸ τοῦ παντὸς ἀρχόμενος συνιστάναι σῶμα ὁ θεὸς ἐποίει. |

---

[5] Respectively, 49a6, 88d6, 51a4.

[6] A 'difficult and obscure passage' concerning a 'mysterious "third kind"' (Zeyl (2010), p. 118).

[7] For instance, cf. Vlastos (1939), p. 71.

[8] Cf. the analysis in Miller (2003), chapter I, in particular p. 32.

[9] Cf. for instance, *infra*, note 52. They may even come from Plato's himself in his lectures. As a matter of fact, one of the strongest Aristotelian criticisms against Plato's theory of Forms is the so called 'third man argument' (*Metaphysics*, I, 9, 990a33-b22, and again in VII, 13, 1039a2-3) which appears already in Plato's *Parmenides* (132a6-b2). In his analysis of three works of Plato, the *Republic*, the *Statesman* and the *Laws* vs. Aristotle's criticisms, M. Vegetti concludes that Plato's works were discussed in the Academy, and Plato answered the objections raised there in his ulterior dialogues (Vegetti (2000), p. 450-1). However, it is reasonable to consider that the same harsh discussions occurred before their completion, and Plato answered them in his works. Of course, some people were not convinced, Aristotle among others, so that the latter used them in his ulterior criticisms.

[10] 'περὶ φύσεως τοῦ παντὸς' (27a4).

[11] 27a6.

[12] 31b4-8, 32b3-4.

[13] The common opinion at that time; cf. for instance G. Morrow: 'Plato, like Aristotle after him, accepts in principle the *current view* that the sensible things of our experience are composed of four kinds of bodies — fire, air, water, and earth — mixed with one another in varying proportions.' (Morrow (1968), p. 13; our emphasizes).

[14] Cf. 48d, where Timaeus warn his audience that they can expect from him 'a strange and unusual narration' ('ἀτόπου καὶ ἀήθους διηγήσεως', 48d5-6), repeated again in 53c1.

| | |
|---|---|
| [synistasthai] just two things in a beautiful way [kalos], without a third; there has to be some bond [desmos] between the two that unites them. Now the most beautiful [kallistos] bond [desmos] is the one that makes this third term and the two others it connects [syndoumena] a unity in the fullest sense; | δύο δὲ μόνω καλῶς συνίστασθαι τρίτου χωρὶς [31c] οὐ δυνατόν: δεσμὸν γὰρ ἐν μέσῳ δεῖ τινα ἀμφοῖν συναγωγὸν γίγνεσθαι. δεσμῶν δὲ κάλλιστος ὃς ἂν αὑτὸν καὶ τὰ συνδούμενα ὅτι μάλιστα ἓν ποιῇ, |
| {and this is accomplished in the most beautiful way [kallista] by proportion [analogia]. For whenever of three numbers, either integers or powers, the middle one between any two of them is such that what the first is to it, it is to the last, and, on the other hand [authis] in turn [palin], what the last one is to the middle, it is to the first. Then, since the middle one turns out to be both the first and the last, and the last and the first likewise both turn out to be the middle, hence, by necessity [ex anagkês], they will all turn out to be in an identical relationship to each other; and, given this relationship, they will form all a single unit.} | {τοῦτο δὲ πέφυκεν ἀναλογία κάλλιστα ἀποτελεῖν. ὁπόταν γὰρ ἀριθμῶν τριῶν εἴτε ὄγκων [32a] εἴτε δυνάμεων ὡντινωνοῦν ᾖ τὸ μέσον, ὅτιπερ τὸ πρῶτον πρὸς αὐτό, τοῦτο αὐτὸ πρὸς τὸ ἔσχατον, καὶ πάλιν αὖθις, ὅτι τὸ ἔσχατον πρὸς τὸ μέσον, τὸ μέσον πρὸς τὸ πρῶτον, τότε τὸ μέσον μὲν πρῶτον καὶ ἔσχατον γιγνόμενον, τὸ δ᾽ ἔσχατον καὶ τὸ πρῶτον αὖ μέσα ἀμφότερα, πάνθ᾽ οὕτως ἐξ ἀνάγκης τὰ αὐτὰ εἶναι συμβήσεται, τὰ αὐτὰ δὲ γενόμενα ἀλλήλοις ἓν πάντα ἔσται.} |
| That being said, if the body of the universe were to have come to be as a two-dimensional plane, without any depth, a single middle term [mesotês] would have sufficed to bind together [syndein] this middle term and the two conjoining terms.<br>As it was, however, the world was to be a solid, and solids are never [oudepote] united [synarmottousin] by just one middle term [mesotês] but always [aei] by two. Hence, the god set water and air between fire and earth, and made them in the same ratio [ana ton auton logos] to one another, as far as was possible, so that what fire is to air, air is to water, and what air is to water, water is to earth. He then bound them together [synedêsen] and combined [synestesato] them to make the visible and tangible universe. As a consequence, to generate the body of the world, these constituents were used, such in kind and four in number, coming into concord [omologêsan] by | εἰ μὲν οὖν ἐπίπεδον μέν, βάθος δὲ μηδὲν ἔχον ἔδει γίγνεσθαι τὸ τοῦ παντὸς σῶμα, μία μεσότης ἂν ἐξήρκει [32b] τά τε μεθ᾽ αὑτῆς συνδεῖν καὶ ἑαυτήν, νῦν δὲ στερεοειδῆ γὰρ αὐτὸν προσῆκεν εἶναι, τὰ δὲ στερεὰ μία μὲν οὐδέποτε, δύο δὲ ἀεὶ μεσότητες συναρμόττουσιν: οὕτω δὴ πυρός τε καὶ γῆς ὕδωρ ἀέρα τε ὁ θεὸς ἐν μέσῳ θείς, καὶ πρὸς ἄλληλα καθ᾽ ὅσον ἦν δυνατὸν ἀνὰ τὸν αὐτὸν λόγον ἀπεργασάμενος, ὅτιπερ πῦρ πρὸς ἀέρα, τοῦτο ἀέρα πρὸς ὕδωρ, καὶ ὅτι ἀὴρ πρὸς ὕδωρ, ὕδωρ πρὸς γῆν, συνέδησεν καὶ συνεστήσατο οὐρανὸν ὁρατὸν καὶ ἁπτόν. καὶ διὰ ταῦτα ἔκ τε δὴ τούτων τοιούτων [32c] καὶ τὸν ἀριθμὸν τεττάρων τὸ |

| means of proportion [*di' analogias*]. Such a proportion bestowed friendship [15] upon it, so that its cohesion is so strong that it could not be undone by anyone but the one who had bound it together. [16] | τοῦ κόσμου σῶμα ἐγεννήθη δι' ἀναλογίας ὁμολογῆσαν, φιλίαν τε ἔσχεν ἐκ τούτων, ὥστε εἰς ταὐτὸν αὐτῷ [17] συνελθὸν ἄλυτον ὑπό του ἄλλου πλὴν ὑπὸ τοῦ συνδήσαντος γενέσθαι. |
|---|---|

There have been many discussions of this passage among modern commentators. One concerns the kind of terms used here ('*arithmos*', '*ogkos*', '*dynamis*' in 31c4-32a1), whether they must be considered as integers, lines or other magnitudes. [18] This, however, does not matter much for the present study, and a detailed discussion of the whole passage is beyond the scope of this article. Let us just say that, more often than not, it is supposed to be considered from an arithmetical point of view, [19] though if it were to be the case it would be hard to see its connection to the rest of a text discussing spatial bodies, which has a clearly geometrical background. [20]

### ii)     The one-middle term proportion

First, Timaeus claims that anything that is born must have a body. Thus, to be corporeal, the universe needs to be both 'tangible and visible'. [21] For the former attribute to hold true, earth is needed, fire for the latter. [22] However, Timaeus claims immediately that the unity of the world needs both elements to be strongly bound, and that the 'most beautiful' [23] bound, [24] is the 'proportion' ('*analogia*') i.e. the equality of ratios between four terms: $a/b = c/d$. [25] This proportion entails the existence of another term, between 'earth' and 'fire', and without further explanation, he remarks that the existence of such a term would not result in a 3-dimensional universe, i.e. with length, width and depth, but a 2-dimensional one. The statement is extremely brief, a clue that it was obvious at Socrates' times. [26]

---

[15] Meaning friendship bonds, perhaps a reference to Empedocles (cf. DK 31 B17).

[16] Thus, the existence of the ordered universe is not absolute, it depends of something outside it.

[17] The same expression appears at 59a6.

[18] Cf. for instance, Festugière (1967), p. 43-44, note 4; Baltzly (2006), introduction, p. 13.

[19] Beginning at least with Proclus (*Commentary on Plato's* Timaeus II, 30.15-33.10), who nevertheless considers a geometrical representation of number *à la* Pythagoras. One has to take into account, as noted by Emilie Kutash, that following Plotinus' emphasis on substantial number, in lines II.30.13-15 'Proclus is giving the priority to number here' over geometry (Kutash (2011), p. 67).

[20] Moreover, in these interpretations, there is no explanation other than some mathematical error by Plato, when he argues in 32b3 that 'solids are never [*oudepote*] united by just one middle term but always [*aei*] by two' (see once again Proclus' *Commentary*, II, 32.15-33.1). Namely, even Timaeus' presentation would be misleading, since instead of the supposed square and cubic integers considered by him, the arithmetical properties he supposedly used, would concern 'similar' integers (ib., 32.1-13).

[21] 31b4-5. A common opinion, since a little further, speaking about the corporeal nature of the elements, he claims: 'First of all, it is plain to everyone, I'm sure …' ('δῆλόν που καὶ παντί…', 53c4).

[22] For light results from fire, and without light nothing is visible (45b-e).

[23] The usual meaning of '*kallistos*', though many scholars, including Zeyl, translate it mostly as 'the best'. However, in view of the importance of symmetry, unity, simplicity and beauty as main characteristics of the universe (cf. *supra*, §3.i), it needs to retain here its usual meaning, especially since 'all that is good is beautiful' (87c4-5).

[24] 31c2-3.

[25] 31c3-4. It was already emphasized in *Gorgias* that 'geometric equality' ('ἡ ἰσότης ἡ γεωμετρικὴ'), i.e. proportion is the bond between humans as well as the whole parts of the universe and even gods (508a). Here, Timaeus does not consider 'proportion' ('*analogia*') in general, but a particular case, the geometric mean i.e. the following equality of ratios between three terms: $a/c = c/b$ ('For whenever of three numbers…', 31c4).

The construction of a geometric mean $c$ between two magnitudes $a$ and $b$ results from Pythagoras' theorem (more exactly, one evident corollary thereof):

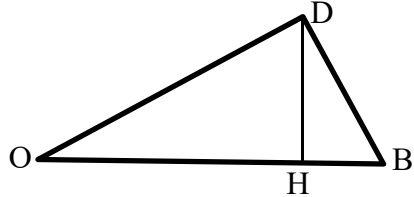

*Figure 1*

In the above drawing, the corollary of Pythagoras' theorem asserts that the height DH = $c$ is the geometric mean between OH = $a$ and HB = $b$. [27] Thus, $a$ and $b$ being given, it is possible (and easy) to obtain geometrically their geometric mean. [28] The importance of this result is that, for any two given lengths $a$ and $b$, it entails the existence of a geometrical mean. [29] In the *Gorgias*, the same property of the 'geometrical equality' [30] is claimed, though without any mathematical explanation.

The most relevant part of this passage for our study is the text within brackets (31c3-32a7). Let us consider it in more detail, using modern mathematical notation for the convenience of the reader. Let $a$, $b$, $c$ be a three term proportion i.e. such that $a/c = c/b$. [31] Though not explicitly stated, an important mathematical result is used: any proportion remains a proportion when the ratios are inversed. The binding of all three terms follows not only from the fact that the ratios $a$ to $c$ and $c$ to $b$ are the same, but from the symmetry of the operation; the possibility of permuting the position of $c$ (and simultaneously $a$ and $b$) and maintaining a proportion. [32] As in the *Gorgias*, it is described as the bond of friendship ('*philia*', 32c2), so strong that it is impossible to be undone, except by a god. The brevity of Timaeus' exposition, which represents a difficulty for us, as it already did for the commentators of late Antiquity, results from the fact that, unlike us, his as well as Plato's audience was 'well schooled' in these matters. [33]

---

[26] And not only at Plato's time as many scholars argue. We have tried to show in previous studies, that it was impossible for Plato, in a city like Athens, where almost everyone knew each other, with so many different highly combative schools and strong thinkers, to describe in writing unlikely situations, in particular wholly anachronistic ones.

[27] In modern notation: $a/c = c/b \Leftrightarrow a \times b = c^2$.

[28] 'What is squaring? The construction of an equilateral rectangle [a square] equal to a given non-square rectangle. Such a definition (…) tells us that squaring is the discovery of a line which is a mean proportional between the two unequal sides of the given rectangle' (Aristotle, *De Anima* II, 2, 413a17-20; J. Smith translation, with slight modification).

[29] This result is described in the *Epinomis* as obviously 'a miracle of divine, not human origin (…) to anyone who can understand it.' (990d4-6; R. McKirahan Jr.'s translation).

[30] Cf. *supra*, note 25.

[31] Integers, lines or other magnitudes, cf. *supra*, note 18.

[32] In modern notation, we have the following equalities: $a/c = c/b$ and $b/c = c/a$. In the second one, $c$ takes the place of $b$ and $a$, while $b$ and $a$ take the place of $c$. The fundamental point is that there is no privileged position for any of these terms, so that there will be no tendency for them to diverge. Symmetry appears here as the most important feature of a proportion.

[33] Cf. *infra*, note 53.

### iii)   The two-middle term proportion

Then, without elaborating, Timaeus claims this solution is correct for planar things, adding that in a plane universe, one geometric mean would be enough, so that only a third term between earth and fire would be needed. [34] However, for the spatial bodies that make up the universe, the problem is different: instead of one middle term, a two-middle term is required. [35] Timaeus does not bother to add any argument to this cryptic statement. Indeed, we, and *a fortiori* a reader at the time of Plato, could hardly avoid connecting it to the famous problem of the 'doubling of the cube'. [36] This problem consists in finding a cube (whose volume is) equal to (the volume of) a given parallelepiped. [37] It is incontrovertible that at the time of Socrates, it was well known that the solution of the previous problem presupposes the demonstration of another one: finding two geometric means between two given magnitudes, a result assigned to the mathematician Hippocrates of Chios. [38] The dating and the output of Archytas' solution is open to discussion. [39] However, rightly understood, this passage gives some essentials hints about this question.

Thus, what does Timaeus say on proportions? First, he claims it is the most beautiful bond between its elements; then he adds that in plane geometry one 'middle term' is enough to bind two terms i.e. this bond is entailed by the existence of a geometric mean between these two terms. [40] Now, what does it mean to claim that when depth is added to length and width, [41] then one middle term is not enough, but two are needed? The answer is straightforward. Between the two terms $a$ and $b$, we now need two geometrical means, $c$ and $d$, such that the two following equalities are verified: $c^2 = a \times d$ and $d^2 = c \times b$; or closer to the language of ancient Greek geometers, the square on $c$ (respectively on $d$) is equal to the rectangle on $a$ and $d$ (respectively on $c$ and $b$). [42] From Hippocrates' result, this provides the solution to the doubling of the cube, [43] or more generally the problem of obtaining a cube equal to any given volume. Now, Timaeus claims two things: the solution of the previous problem is impossible in a plane, [44] and spatial bodies must be used, as done in Archytas' solution. Hence, according

---

[34] 32a7-b2.

[35] 'As it was, however, the universe was to be a solid, and solids are never joined together by just one middle term but **always** by two.' (32b2-3). In modern notation, instead of a proportion of the form: $a/c = c/b$, we need a two-proportion set: $a/c = c/d = d/b$.

[36] A connection commonly made in Antiquity: cf. the Democritus (3rd century CE) as quoted by Proclus in his *Commentary of the Timaeus* (II, 33.15-20).

[37] In modern notation, if the sides of the parallelepiped are respectively $a, b, c$, find geometrically a cube of side $d$ such that $d^3 = a \times b \times c$. The problem of 'doubling' consists in finding a cube twice (the volume of) a given one; however, Archytas' solution solves the general problem (though Archytas' authorship is not uncontroversial, it is traditionally assigned to him, and we will refer to it as Archytas' solution).

[38] It is used in Archytas' solution, and likely obtained no later than the mid-5th century BCE. A reasonable dramatic dating of the dialogue is 430-420 BCE. Debra Nails proposes the more precise dating of 429 BCE (Nails (2002), p. 324-326). For the mathematical details, we refer to the **Annex I** in Brisson-Ofman (to appear)

[39] Cf. for instance, Brisson (2013)

[40] Cf. *supra*, note 27.

[41] In other words, when we consider spatial figures instead of plane ones.

[42] This results from the two proportions: $a/c = c/d$ and $c/d = d/b$. Moreover, the proportion $a/c = d/b$ entails that $a \times b = d \times c$; thus $c^2 = a \times d$ entails $c^3 = (a \times d) \times c = a \times (d \times c) = a \times (a \times b) = a^2 \times b$. Thus, for $b = 2a$, we get: $c^3 = 2a^3$.

[43] According to the last equality of the previous note, the cube of side $c$ solves the problem of doubling the cube of side $a$.

[44] At least when the usual tools of a geometer, the rule and the compass, are used. Other mechanical tools were quickly invented by Greek mathematicians to obtain practical solutions of the problem (cf. Eutocius'

to our argumentation, [45] this solution was well known, not only at Plato's time, but at the dramatic dating of the dialogue, i.e. around 430-420 BCE.

Thus, continues Timaeus, between earth and fire, two other elements are needed: water and air. Since the universe is not a mathematical one, and does not follow a strict rational order, the three ratios [46] may not be exactly the same but as close as possible. [47] Once again, no explanation is given for the use of these two new elements, [48] which implies it was the common opinion in the 5$^{th}$-4$^{th}$ centuries BCE. As we can see, Timaeus' purpose is to give a mathematical framework in which the usual 'physics' of his time takes place. [49]

Here (32c4) ends the first presentation of the mathematical structure of the universe, by means of the then common four-element framework. Further, Timaeus says it is needed to clarify the above construction and explain their dynamic, as the four elements appear to be morphing into each other in a continuous circle, from birth to destruction and conversely. [50] However, in a second part, he adds that this view is too simple, and needs to be corrected. [51] Thus the supposed general circle of birth and destruction of the four elements needs to be modified, since it is inconsistent with the mathematical principles underlying the general construction of the universe. [52]

## 3. A two-triangle universe

In the passage beginning at 49b1, after a short introduction on the difficulty of the task, Timaeus considers the four elements in more detail (53c4-60c2). He notes that his

---

"Commentary on Archimedes' On the Sphere and Cylinder", II, in Heiberg (1881); cf. also Huffman (2005), III, 3, p. 342-401). This claim is all the more amazing that, while it is true, its proof was not obtained till the 19th century, through some results connected to Galois' theory. It shows the deep geometrical understanding of ancient Greek geometers. Indeed, a solution by 'rule and compass' was not only a question of tools, and the ancient Greek geometers certainly did not think in this way. They considered two elementary geometrical figures, the line and the circle, and all geometrical figures were compounded of them, as claimed by Parmenides in Plato's eponymous dialogue (137d8-e1, 145b3-5). This point is also emphasized by Aristotle: 'all movement that is in place, all locomotion, as we term it, is either straight or circular or a combination of these two, which are the only simple movements. And the reason of this is that these two, the straight and the circular line, are the only simple magnitudes. ('Πᾶσα δὲ κίνησις ὅση κατὰ τόπον, ἣν καλοῦμεν φοράν, ἢ εὐθεῖα ἢ κύκλῳ ἢ ἐκ τούτων μικτή', *On the Heavens* I, 2, 268b17-18)'. It was therefore possible, at least theoretically, to construct any figure by rule and compass.
[45] Cf. *supra*, note 26.
[46] In modern notation, *a/c*, *c/d*, *d/b*.
[47] 32b3-5. This follows from the existence of 'necessity', which even the gods need to take into account.
[48] Except much further, where it is stated that these four elements are the most beautiful bodies (53e5), cf. *infra*, note 77.
[49] See 68d2-7, and also *infra*, note 163, for the mistake of thinking that what we see matters more than the fundamental (mathematical) principles; cf. also Cornford (1937), p. 45. This is consistent with Simplicius' quotation from Sosigenes, to the effect that in astronomy, Plato was always looking for a theory that would 'save the phenomena' ('σώζειν τὰ φαινόμενα') (*Commentary on Aristotle's* Treatise on the Heavens, II 12, 292b10 (488.3-24)). See also *infra*, note 79.
[50] 49c6-7.
[51] Timaeus' use of the common theory of nature may follow a similar pattern to his approach of the genealogy of the gods. When we lack knowledge about something, the simplest path is to follow the self-proclaimed wise men (those 'who claim to be the offspring of gods', 40d8) i.e. the ancient poets, as ludicrous their claims may be (40d6-e3); namely, one must 'follow custom and believe [the wise]', as long as there is nothing better. Cf. also, *supra*, note 49.
[52] Cf. *infra*, §3.ii). This is paradoxically confirmed by Aristotle's harsh criticisms of Plato or the Platonists (*On the Heavens* III, 7, 306a5-15).

exposition will be different from the usual ones. [53] Indeed, he is breaking not only with all the 'natural thinkers' of the past, but also with almost all physicists in Antiquity and even much later ones, as he proceeds to set the universe on geometrical foundations. Plato's readers are assumed to be familiar with the same fields as Timaeus' 'well-schooled' interlocutors, the different parts of mathematics, particularly geometry. His first statement is that all the elements are bodies with length, width and depth. However, instead of dealing with spatial bodies, he will consider plane surface and the transformation of the elements are based on surfaces, not on volumes. A common error must be avoided: as Cornford asserts, these surfaces are just that, not 'thin plates of corporeal matter, forming boxes with a hollow interior'. [54] This appears so strange to modern commentators that many have added, with many a detail, what happens in Timaeus' universe to volumes or even weights of spatial bodies. Nevertheless, clarifying Timaeus' choice is essential for understanding his presentation of the foundations of the universe. [55]

### i) The beauty of the universe

Timaeus' universe is characterized by *beauty*, and hence by *unity*, *harmony*, and *symmetry*. Right from the outset, Timaeus makes clear that his guiding principle will be beauty. [56] Because the universe must be beautiful, [57] it has to be modeled on the 'what is always changeless'. [58] It is 'of all the things that have come to be (…) the most beautiful'. [59] Again, it possesses a soul ('*psychê*') because it needs intelligence ('*nous*'), for an intelligent being is more beautiful ('*kallion*') than one without intelligence. [60]

For the same reason, all that is visible must be included in the same unique beautiful universe, which must be in a state of order and harmony. [61] Moreover, all living beings have to belong to the same universe 'for anything incomplete could not ever turn out beautiful'. [62] Again, because the universe must be the most beautiful, it must contain all intelligent beings. [63] Hence, beauty, along with order and harmony, are Timaeus' arguments in support for the unicity of the universe. [64] They also point to a strictly limited number of basic components, which explains why he took over the common theory of the four elements, [65] rejecting the atomists' theory of unlimited basic elements. [66] The (geometric) proportion is the most

---

[53] 'My account will be an unusual ('*aêthei*') one, but since you are well schooled in the fields of learning in terms of which I must of necessity proceed with my exposition, I'm sure you'll follow me.' (53c1-3).

[54] Cornford (1937), p. 229; see also Martin (1841).

[55] Aristotle blamed the Platonists because they 'divide bodies into planes and **generate them out** of planes' (II 4, 286b27-28). Cf. also 298b33 and *infra*, notes 94 and 129.

[56] For the modernity of such an approach, cf. *infra*, note 69.

[57] 'καλὸν ἐξ ἀνάγκης' (28a8).

[58] The intelligible world.

[59] 'ὁ μὲν γὰρ κάλλιστος τῶν γεγονότων' (29a5).

[60] 30b1-4.

[61] From a world that was 'not at rest but in discordant ('*plêmmelôs*') and disorderly motion, it was brought (…) to a state of order.' (30a4-5).

[62] 30c5.

[63] 30d3-31a1.

[64] Against the possibility of many or an infinity of universes (31a2-b4).

[65] Cf. also *supra*, note 48.

[66] Cf. Aristotle's statement opposing the atomists, in particular Leucippus, to Plato: 'there is this great difference between the theories [of Plato and Leucippus]: the "indivisibles" of Leukippos (…) are characterized by an infinite variety of figures, while the characterizing figures employed by Plato are limited in number.' (*On Generation and corruption*, 325b27-29). It is not clear why some scholars consider that Aristotle is here bringing both thinkers closer (cf. Vlastos (2005), p. 67).

beautiful bond between its terms because it imposes the most complete and strongest unity between them, making their dissolution impossible. [67] Namely, as we will see in the next paragraph, even the choice of the geometric figures used as fundamental particles on which the whole universe is founded, results from such a search for beauty. [68]

In ancient Greece, beauty is closely connected to unity and symmetry. [69] Symmetry [70] is an essential characteristic of Timaeus' cosmology. In particular, it is because of the symmetry of the sphere 'with its center equidistant from its extremes in all directions', that it is 'absolutely perfect', [71] the figure 'which of all figures is the most perfect and the most self-similar', [72] and always 'the similar is infinitely more beautiful than the dissimilar'. [73] It is the reason for the spherical form of the universe, [74] the reason why its only movement is the uniformly circular one, [75] 'deprived' of all six asymmetrical movements. [76] Later in the dialogue, Timaeus claims that it is again beauty that decides, among all other possibilities, on the choice of the four fundamental elements (earth, water, air, fire). [77] More generally, he adds that there are two kinds of cause. Beautiful and good things are produced by the prime or real causes, [78] while secondary or necessary causes produce haphazardly. [79] Hence, beauty is the best lead to science, thus, the best path to follow to understand Timaeus' account.

Further, summarizing what he has said, he emphasizes the symmetry of the most beautiful possible universe. [80] He claims at length, against the common opinion, that it is impossible to speak universally either of high and low, left and right, or before and behind. [81] As in modern

---

[67] 31c2-4.

[68] This is recalled in several parts of the dialogue in particular at 69b. The main means to achieve beauty is the geometrical symmetry provided by the geometric proportions (cf. also *supra*, note 32). This point is rightly emphasized in Gregory's introduction to his translation of the *Timaeus* (Gregory (2008), p. xv).

[69] Cf. for instance 87c and also *Philebus*, 64d-66b. Let us remark that the importance of symmetry for modern physics cannot be overstated (cf. previous note). For a physicist's point of view, cf. for instance Zee (1986) who claims in his abstract that for 'today's theoretical physicists (…) [beauty] provides us with **an important tool for the exploration of nature**.' (Gross (1996), p. 14258; our emphasis). This point was already noticed by J. Cook Wilson (Wilson (1889), p. 49). Nevertheless, some scholars have considered that this does not justify the added complexity of Timaeus' construction (for instance, Cornford (1937), note 2, p. 217). However, as noted by Bruins, 'if Timaeus divides the square into four triangles through the drawing of their diagonals, and the equilateral triangle into six partial triangles, there a fundamental necessity to do so.' (Bruins (1951), note 2, p. 277); cf. also, Lloyd (2006), p. 462. We will see that this complexity is indeed a necessary property, playing an essential role in the ordered organization of the Platonic cosmos.

[70] The Greek term '*symmetros*' means 'commensurable' in mathematics, but more generally 'in right measure', 'in due proportion', and 'symmetrical'. When Timaeus uses the term outside of a strict mathematical context, especially when he speaks about beauty, he includes the symmetrical aspect, as asymmetry is associated to ugliness ('*aischros*'), disequilibrium ('*paraphorotês*'), disproportion ('*ametros*'), cf. 87e2-5.

[71] 'ἁπάντων τέλεον' (33a7).

[72] 33b6, for the sphere is the most symmetrical geometrical figure.

[73] 'μυρίῳ κάλλιον ὅμοιον ἀνομοίου' (33b7).

[74] 34b1.

[75] 34a3-4.

[76] 34a4-5, i.e. 'forwards and backwards, back and forth, to the right and the left, upwards and downwards' (cf. 43b), the 'wandering' and 'irregular' movements.

[77] 'We should now say which are the most beautiful four bodies that can come to be.' (53d7-e1).

[78] 46e3-4.

[79] The ones that are 'deserted by intelligence, [which] produce only haphazard and disorderly effects every time.' (46e5-6). Nevertheless, they must be understood, if only to enable understanding of the first ones (69a1-5).

[80] A universe where all the things are 'proportionate and symmetrical ('ἀνάλογα καὶ σύμμετρα').' (69b5).

[81] 62c-63e.

physics, symmetry is everywhere in *Timaeus*' cosmology, from the spherical form of the total universe to the imperceptible basic right triangles, as we will see in the next paragraph.[82]

### ii) The right triangles

First Timaeus remarks that any solid is enveloped by surfaces and any rectilinear surface is composed of triangles (53c7-8).[83]

- Then, Timaeus carries out a first reduction: these general (rectilinear) surfaces are considered as addition of triangles.[84]

- The second reduction concerns the triangles themselves: only right triangles are taken into account, for all triangles are the union of two right triangles,[85] since a height divide a triangle in two right triangles.[86]

- The reason for the third reduction is, once again, beauty, that is harmony and symmetry. The goal is to choose the most beautiful among all right triangles. Timaeus divides the former into isosceles and scalene. The isosceles triangle, as the most symmetrical, must be used for the formation of the elements. However, while there is only one isosceles kind of right triangles, there are unlimited kinds of scalene.[87] The distinction between these two kinds of triangles does not constitute a correct division,[88] hence the need to find a model for the scalene triangles, to obtain two balanced genera. It will be the most beautiful among the latter, the half-equilateral triangle.[89] Once again, symmetry is the reason for this choice.

Thus, all four elements derive from two triangles: the isosceles right triangle and the half-equilateral one.[90] Then, however, Timaeus notes the need to correct a previous mistaken claim.[91] Though there seems to be a continuous cycle of transformations between the four elements, this is not true. For there is only one element that arises from the isosceles right triangle, the three others resulting from the scalene. Since it is impossible for an isosceles triangle to become scalene, only the three latter elements may follow a circle of transformations.[92] However, surprisingly, so far it is still not known how the four elements are linked to these two triangles.

---

[82] Plato has made clear 'that symmetry and proportion were introduced down to the smallest detail' (Cornford (1937), p. 239).
[83] Excluding round figures.
[84] In other words, a union of triangles joined by their sides.
[85] 53c8-d2.
[86] This is not true for any height in any triangle, at least when one of its angles is obtuse. However, even in this case, there is always one suitable height, as shown below:

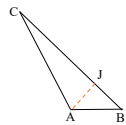

[87] 54a1-2.
[88] Cf. *Statesman*, 262a-263a, and also Brisson-Ofman (2017), p. 129-132.
[89] 54a2-3.
[90] 54b1-5.
[91] 'First, we see the thing that we have just now been calling water condensing and turning to stones and earth. Next, we see this same thing (…) turning to wind and air, and air (…) turning to fire. And then we see fire (…) turning back to the form of air, and air (…) turning back into cloud and mist (…) [then] we see them turning into flowing water, which we see turning to earth and stones once again. In this way, then, they transmit their coming to be one to the other in a circle, as we see.' (49b7-c7). In the passage, earth has a special status, since contrary to the other elements, it is not said it is turning back into another one i.e. water. Thus the cycle is left incomplete, perhaps suggesting already no one ever has ever seen such a phenomenon, thus opening the way to a later correction (cf. *infra*, note 103).

### iii) A world of regular polyhedrons

Since the elements are entirely defined by triangles, it remains to explain how these two right triangles give rise to the four elements. Going further on the path of the mathematisation of the universe, Timaeus presents a one-to-one correspondence between the four elements and four regular polyhedrons,[93] though still without giving the correspondence between the polyhedrons and each element. These polyhedrons are differentiated by two characteristics of their faces:

- Their number: 6 squares for the cube, 4 equilateral triangles for a tetrahedron, 8 for an octahedron and 20 for an icosahedron.
- Their form: squares for a cube and equilateral triangles for the others.

A square is formed by the addition of two isosceles right triangles, an equilateral triangle by two half-equilateral right triangles. Hence, the correspondence between the elements and the triangles is given by means of the formation of the faces of the regular polyhedrons by these triangles.[94] Accordingly, the element corresponding to the cube will not be transformed into another element, contrary to the three other elements, whose faces are equilateral.[95] Hence, one would expect that two right triangles (either both isosceles or both scalene) will form each of the faces of the polyhedrons, thus of the four elements.[96]

## 4. The two-triangle construction

The concern for unity, which appears at the very beginning of Timaeus' account,[97] leads to the choice of the right triangles as foundational figures. Therefore, a natural construction would be to consider the elementary right triangles as half-square and half-equilateral. In both case, the complete figure (square and equilateral triangles) would be the result of the addition of two triangles (respectively the isosceles and the scalene), as seen below.

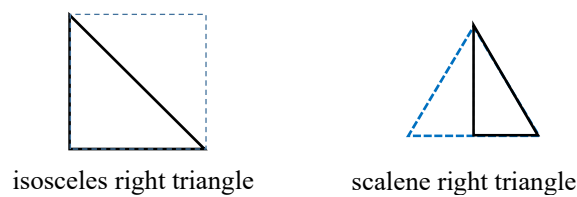

isosceles right triangle     scalene right triangle

Figure 2

---

[92] 54b8-c5.
[93] The fifth, the dodecahedron, is excluded from the list, but it will appear later.
[94] It is essential, though often forgotten, that the elements are exclusively characterized by their faces, not their volumes. Cf. *supra*, §3 and note 55, and also *infra*, note 129.
[95] For a beautiful representation in perspective of the five regular polyhedrons by Leonard de Vinci, cf. Pacioli (1509), plates ii, viii, xviii, xxii, xxviii.
[96] Cf. *infra*, Figure 2.
[97] Cf. 31c2-3 and *supra*, note 24. Cf. also *infra*, note 100.

However, Timaeus does not follow this way. He uses four isosceles right triangles to obtain the squares, i.e. the faces of the cube, and 6 isosceles right triangles for the equilateral triangles, i.e. the faces of three other regular polyhedrons. [98] Hence, a question that has puzzled the commentators, from the Antiquity to today: the apparent unnecessary complexity of Timaeus' geometrical construction. [99] Why does he choose four isosceles and six scalene right triangles, while two are enough in both cases? The path to a solution is to focus on the main character of Timaeus' universe, that is, beauty, unity and symmetry. [100]

Let us compare the two geometrical constructions, i.e. the simplest and more 'economical' vs. that of Timaeus (cf. figure below):

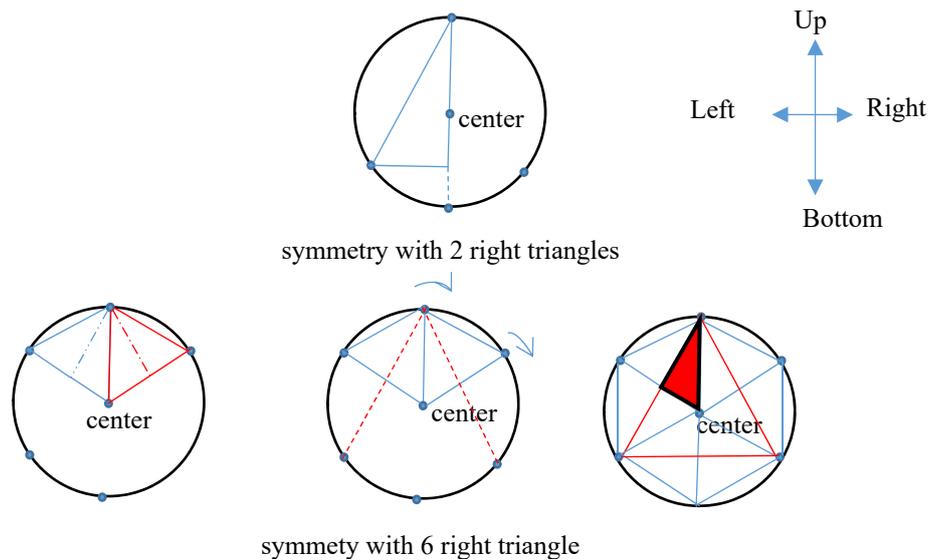

Figure 3

The difference in symmetry of the figures is obvious. In the first drawing, the most economical, [101] there is one axis of symmetry, therefore one line and one direction are privileged. In Timaeus' construction, the symmetry is identical with respect to the three main lines that derive from the equilateral triangle. The figure consists in the rotation on the circle of two smaller equilateral triangles, as seen on the second, lower drawing. Here, we obtain a complete circular symmetry: Timaeus' six scalene right triangles are obtained from the intersection of the equilateral triangle and the six diagonals of the hexagon inscribed in the same circle. This is obviously not the case for the more economical upper construction. [102]

The same opposition appears for the isosceles triangles forming a square. By using two triangles to obtain a square, one of the diagonals is privileged, while in Timaeus' construction

---

[98] 54d-e.

[99] 'Here we encounter one of the points that we noted as never having been satisfactorily explained.' (Cornford (1937), p. 234); 'something of a mystery' (Zeyl (2000), note 14, p. 34); one of the two 'major puzzles which Cornford hopes to solve', but failed (Pohle (1971), p. 38); a 'mystery' for Lloyd (2006) (p. 4).

[100] Cf. 30a and *supra*, §3.i), in particular notes 68 and 69.

[101] In line with Pohle (1971), this construction will be called in the remainder of the article, the 'economical construction' as opposed to Timaeus' one.

[102] Cf. *supra*, §3.i). Francis Cornford proposes another explanation for Plato's use of these right triangles rather than directly equilateral triangles or squares (cf. *infra*, §5.i).b)). However, his construction does not follow that of Timaeus, but is based on propositions IV.15, XIII.12 and XIII.13 of Euclid's *Elements* (Cornford (1937), p. 235-238). Yet, the figure illustrating his argumentation highlights one fundamental point, the symmetry of the construction (p. 236). It is indeed essentially the last drawing of Figure 3 above.

both have the same role, and the symmetry is total.

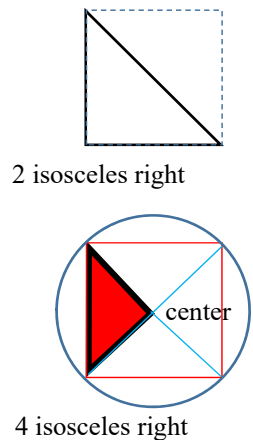

2 isosceles right

4 isosceles right

Figure 4

### 5. Unlimited diversity vs. a two-triangle universe

Only after providing the construction of the polyhedrons from the elementary surfaces (squares and equilateral triangles), is Timaeus to give the correspondence between each of the four elements and the regular polyhedrons: tetrahedrons are components of fire, cubes of earth, octahedrons of air, and icosahedrons of water. It is geometrically impossible for the earth to be transformed into any other element, unlike the three others.[103] Timaeus gives some examples of the laws of transformation. A particle of water may be decomposed into one particle of fire and two of air, for the faces of the first, as an icosahedron, are composed of 20 equilateral triangles; the fire, as a tetrahedron, of 4 such triangles; and the air of 8 of them.[104] Another possibility is that one particle of air (8 equilateral triangles) will give two particles of fire (4 such triangles each),[105] or, conversely, two particles of fire will give one particle of air.[106] In the last example, from 2 + ½ ('δυοῖν ὅλοιν καὶ ἡμίσεος') particles of air, one will obtain one particle of water.[107]

But a daunting problem now arises: how, from only **two** right triangles, or even four polyhedrons, is it possible to get the extraordinary diversity of the universe? An answer is given in a further passage,[108] where Timaeus explains the 'unlimited' diversity of sensible things. It results from, on the one hand, the diversity of the sizes of the particles belonging to the same element, on the other hand, the combination of the elementary particles. [109] Until this point, Timaeus considered only four primitive elements (earth, water, air, fire) associated one-to-one to four regular polyhedrons (cube, icosahedron, octahedron, tetrahedron). [110] Now, as he gets into more details, he claims that, within each elements, there are many kinds of different sizes.

---

[103] 'the parts of earth will never pass into another form.' ('γῆ γένοιτο—οὐ γὰρ εἰς ἄλλο γε εἶδος ἔλθοι ποτ' ἄν', 56d5-6). Thus, it is only now that one understands the problem of the mistaken cycle between the four elements (cf. *supra*, notes 91 and 92). While water, air and fire can be turn into each other, earth will always be earth; it cannot change into any other element. The absence of the earth from the cycle of elementary transformations was certainly not a common view, cf. for instance, cf. Aristotle, *Generation and Corruption* (II, 4, 331a12-20), *On the Heavens* (III, 2, 306a2-9). Moreover, Timaeus himself previously said: 'we see (…) water turning to stones and earth, (…) and flowing water turning again to earth' (cf. note 91).

[104] 20 = 1 × 4 + 2 × 8 (56d7-8). Of course, there are many other possibilities even simpler, for instance, the same particle of water could be decomposed into 5 particles of fire (20 = 5 × 4).

[105] 56e1-2.

[106] 56e2-5.

[107] (2 × 8 + ½ × 8) = 16 + 4 = 20 equilateral triangles (56e6-7). The use of fractions of particle is another testimony that the laws of transformation are not essentially about the particles themselves but the triangles i.e. their faces.

[108] 58c5-60c4.

This statement seems clearly at odds with Timaeus' former claim that the different elements can transform into each other. This is corrected later, when he adds that only those corresponding to the polyhedrons whose faces are equilateral triangles, i.e. excluding the element earth, participate to this cycle of transformations. [111] Timaeus goes so far as to give explicit examples of such intra-elementary transformation. That said, if the size of the elementary particles is not the same for all the three elements, it is hard to see how such transformations could take place: particles of any elements will not be decomposed into any particle of another element, but only into particles of the same size, in contradiction with the circle of transformation. [112] Even worse, diversity in the same element, including this time the element earth, would also be impossible if the different kinds of particles within the same element (for instance, copper and ice, which both belong to the element water) are of different sizes. [113] Moreover, the elementary triangles (isosceles or scalene) would not combine with each other, but would need to 'wait', so to speak, until they meet triangles of the same size.

One does not, however, find any comment in Aristotle against such blatant inconsistencies, whereas what seem to be much less obvious mathematical problems were harshly criticized by him. [114] Likewise, one does not find either any defense of these supposed difficulties by Plato's later supporters, whereas, once again, much pettier ones were addressed at length. [115] Hence, this was not considered an issue during Antiquity. It is a clue not only that this was not considered as a hindrance to Timaeus' construction, but that the 'solution' [116] was obvious to a reader of this period.

In order to clear any ambiguity, one will call in what follows 'basic' right triangles, the two right triangles (one half-equilateral, the other isosceles), making up the 'basic' faces (the 'basic' equilateral triangles or squares) of the four (regular) 'basic' polyhedrons associated to the four 'basic' elements; and 'elementary' triangles, squares or polyhedrons, the triangles, squares or polyhedrons of different sizes i.e. the faces of the larger or smaller particles of each element. [117]

---

[109] For the diversity of varieties of fire, cf. 58c5-d1; of air, 58d1-4; of water 58d4-60b5; of earth 60b6-61c2.

[110] Cf. *supra*, §3.iii). Only then, one learns the assignment of the regular polyhedrons to the elements: 'Let us now 'assign' ('*dianeimômen*') to fire, earth, water and air the structures which have just been given their formations in our speech' (55d6-8).

[111] Cf. *supra*, §3.ii), in particular note 89.

[112] Cf. *supra*, note 91.

[113] Namely, in this case, as Cornford states 'there will thus be (…) parallel processes of transformations' ((1937), p. 232). Nothing on that can be found in the Timaeus.

[114] For instance, the impossibility for regular polyhedrons to fill completely a space, except, says Aristotle, for cubes and tetrahedrons (namely, only cubes are able to do this). We will return to this problem in a subsequent article.

[115] For instance, see Proclus' *Commentary on Plato's* Timaeus.

[116] Of what appears as a problem, but only to us, modern readers.

[117] 'Each of these two constructions did not originally yield a triangle that had just one size, but triangles that were both smaller and larger ('ἐλάττω τε καὶ μείζω')' (57d1-2).

### i) Cornford's solution

These difficulties are well highlighted by Cornford,[118] who proposes an elegant solution.

a) The size of the elementary triangles can vary because they are able to join side by side. Cornford remarks that Timaeus' construction, contrary to the 'economical'[119] one, allows adding another one to this sequence, for two identical isosceles right triangles form another isosceles right triangle whose side is multiplied not by 2 but by √2, as shown in the figure below.

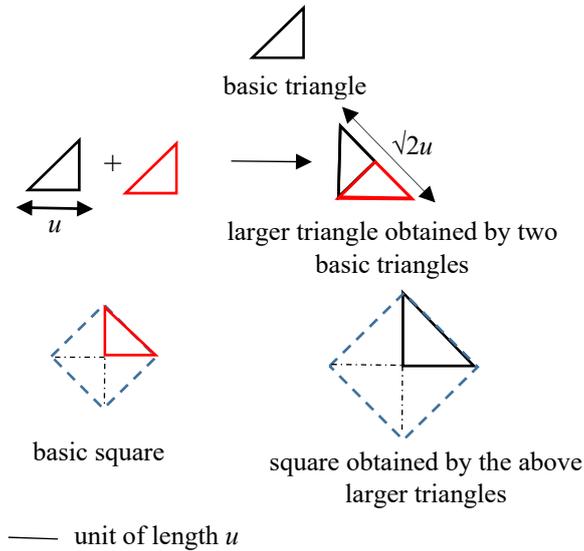

Figure 5

In the second row of the above figure, the ratio of the sides of the larger (similar) triangle to the sides of the basic triangle is √2; hence, there is the same ratio between the sides of the two squares in the last row. Thus, combining the two series of squares, one gets the figure below:

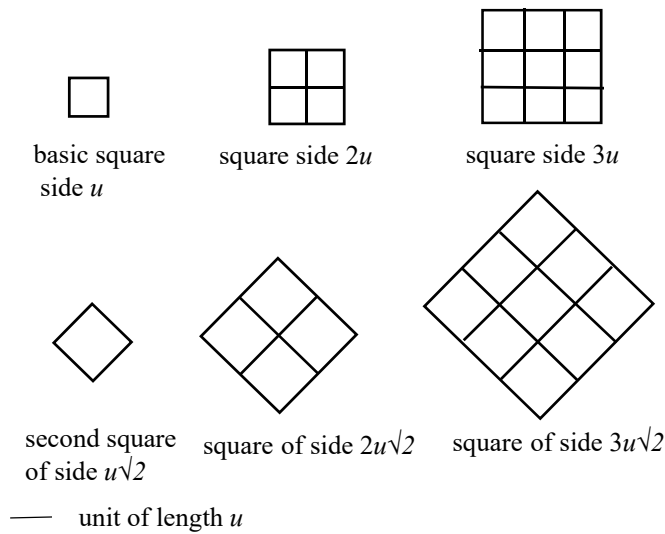

Figure 6

---

[118] Cornford (1937), p. 231-233.
[119] Cf. *supra*, note 101

Thus, in the above figure, one obtains a combination of the sequence of squares (faces of the cube) of side $nu$ and $nu\sqrt{2}$. Hence, in term of surfaces, one has a sequence of squares of areas $n^2$ and $2n^2$ times the basic square ($n$ integer and $u$ side of the basic square).

For the scalene triangles forming the equilateral triangles, faces of the other three regular polyhedrons, it is possible to proceed in a similar way, using three such basic scalene triangles. Thus, one obtains a new (similar) triangle whose sides are multiplied by $\sqrt{3}$ and the area by 3, [120] as shown in the figure below. [121]

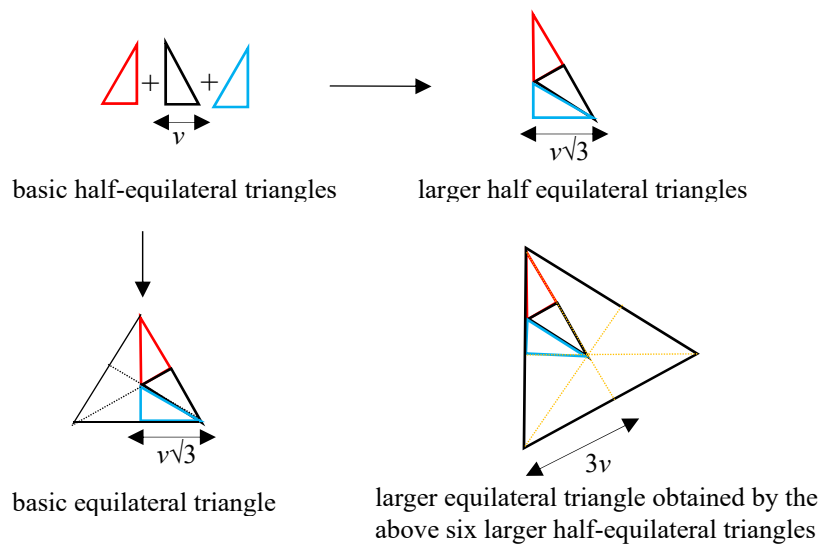

basic half-equilateral triangles      larger half equilateral triangles basic equilateral triangle      larger equilateral triangle obtained by the above six larger half-equilateral triangles

Figure 7

Then, one gets a sequence of triangles of sides $nv$ and $nv\sqrt{3}$ and of areas $n^2v^2$ and $3n^2v^2$ ($n$ integer and $v$ the smallest side of the basic half-equilateral triangle).

Cornford claims then that his solution explains both why there is no inconsistency between, on the one hand, the diversity of the elements and, on the other hand, the necessity of only one basic isosceles right triangle and one basic scalene right triangle (consequence of the cycle of transformation), the reason of Plato's choice of Timaeus' construction over the economical one, being the former gives more kinds of each elements than the latter. Though some leading scholars had endorsed it, sometimes with some reluctance, [122] it is not without problems, so that many others were not convinced. [123] Thus, let us consider some of the difficulties underlying it.

  b)     Let us first begin about his assertion that his solution explains Plato's choice (four isosceles and six scalene right triangles) rather than the economical one (two such triangles). Since Cornford seems sometimes to hesitate between opposite argumentations, we will consider them successively.

    1) He first claims that for each element there are few different kinds ('grades') issued from a basic one. [124] However, this solution would not explain the 'endless diversity' of

---

[120] The economical construction would only lead to triangles of sides $nv$ and of areas $n^2v^2$.
[121] Thus, one gets new larger particles of the three other elements, fire, air and water.
[122] Cf. for some references, Pohle (1971), p. 36-37; Lloyd (2006), p. 4.
[123] For instance Pohle (1971), in particular p. 43-44.
[124] Cornford (1937), p. 231.

natural things. [125] To avoid this difficulty, instead of beginning with the basic triangle, he introduces an arbitrary isosceles right triangle, and shows that it is easily divided in two similar triangles as many times as is required. [126] Since he agrees that there are basic right triangles, of minimal sizes, [127] then again, the diversity of natural things would be severely limited. [128]

2) A little further, the author re-addresses the question, and recognizes that according to the economical construction, his solution leads to elementary particles that would be too large to be invisible. [129] However, he claims that using Timaeus' construction instead of the economical one, his interpretation avoids this difficulty. Unfortunately, as seen above, the increase in capacity of the growth within Cornford's solution is not great: between one basic elementary particle and any kind of elementary particle of given size, there would be less than twice possible kinds of elementary particles, according to Timaeus' construction with respect to the economical one. [130] Thus, Cornford's claim that Timaeus' construction is better adapted for the diversity of natural things is questionable. Namely, since according to him, there is a limited number (either small or large, it does not matter) of different kinds for each element, let $u$ and $v$ be as in the two above figures. As Cornford remarks, all the particles of the different kinds for each element must be invisible. Then, changing $u$ for $u/2$ and $v$ for $v/2$, all elementary particles in the economical construction would be smaller than in Cornford's solution. Thus, if, as he claims, the latter construction leads to elementary particles too large to be invisible, the same would apply for his own solution. [131]

3) Another problem is Cornford's understanding of the passage 57d. He claims that the 'endless' diversity of Nature [132] is not consequence of an unlimited kinds of elements that are, according to him, in limited number, but of the 'indefinitely numerous combinations of all these varieties in composite bodies' i.e. of the mixing of these different kinds in spatial bodies. [133] However, nothing in the text suggest such dissymmetry between a limited number of kinds of elements and an 'endless' number of compositions of elements. Namely, if that was the case, the existence of 'smaller and larger' particles would be useless, and one kind of each element would explain the 'infinite' diversity as well, even though. And yet, to explain the diversity of the nature, Timaeus puts first the multiplicity of the necessity of different kinds within each element. [134] Hence, it is more reasonable to understand that Timaeus claims that the kinds of each element are 'endless' as well.[135] Thus, let us consider what happens for Cornford's solution in this case. Namely, explaining probably Cornford's rejection, the consequence is even worse: how small the basic triangles or squares may be, it leads to some particles larger than any fixed number, even larger than the whole universe.[136]

---

[125] Ibid, p. 230

[126] Ibid, p. 233. His solution for scalene triangles is more complicated, but based on the same idea.

[127] This is besides obvious if one hypothesizes a finite number of kinds of the four elements (ibid, p. 231).

[128] Cf. *infra*, note 134

[129] *Ibid.*

[130] Cf. *supra*, a).

[131] Besides, since human acuity is so variable, it would be strange that, by doubling the size, from invisible the particles become visible.

[132] Cf. *supra*, note 125.

[133] *Ibid.*, p. 231.

[134] It could be argued that if the number of kinds of the elements is finite, even if it is a very large number, it is always possible to start from basic figures (squares or equilateral triangles) so that all the particles are extremely small, thus invisible (see for instance Artmann-Schäffer (1993), p. 258). However, since then the size of the basic particles would depend of this number, it (or at least some upper-bound) would need to be known before the construction of the universe. Such a strong hypothesis severly limiting of the possibilities of the diversity of the nature, is completely lacking in Timaeus' account, and is inconsistent with his description of the universe. Namely, Cornford rightly does not support it for he notes that, at least in the economical construction, the juxtaposition of triangles will necessarily lead to discard the invisibility of the elementary particles (ib., p. 238).

Thus, in both the finite and the 'endless' cases, it is impossible for the elementary particles to remain invisible, in open contradiction with Timaeus' account.

c)  Let us now consider a weaker hypothesis: Cornford's solution as only an answer to the problem of the diversity of the universe, independently of the problem of Timaeus' construction vs. the economical one. Yet even this weaker assumption is problematic.

1) The sizes of the particles of the different kinds of the four elements are respectively of the form $nu$ and $nu\sqrt{2}$ for the earth, $nv$ and $nv\sqrt{3}$ for the third other elements ($n$ integer).[137] Thus, Cornford's solution entails a new series of ratios between the elementary particles,[138] which is nowhere to be found in Timaeus' account.

2) The sequence of elementary squares has two generators (squares of sides $u$ and $u\sqrt{2}$) instead of one,[139] and the same applies for the equilateral triangles (equilateral triangles of sides $u$ and $u\sqrt{3}$).[140] Hence, the construction loses its unity and symmetry, contrary to Cornford's own claim about the fundamental character of symmetry in the *Timaeus*.[141]

3) However, one of the two main difficulties is that the addition of a second series does not change much with regard to the problem of growth, for it multiplies the number of different squares or equilateral triangles by a number less than two.[142] Since its purpose is to solve the problem of the 'unlimited' diversity of things in nature, [143] the number of sizes of these figures must be very large, even 'unlimited', so that the areas of the different elementary polyhedrons would become quite gigantic. It would soon become impossible for the elementary particles to remain invisible, [144] a difficulty already met by the atomists. [145]

---

[135] It was the classical interpretation before Cornford's analysis (for instance, cf. Martin (1841), note lxxvi, p. 254; Rivaud's French translation of the *Timaeus* (1925); see also Mugler (1960), p. 22-26)). It does not mean that there are infinite different kinds of each element, but that there is no fixed limit (cf. *infra*, note 156), as needed according to Cornford's reading (cf. also *supra*, 1)). G. Vlastos ((1967), note 8, p. 205) devotes a long note on the meaning of 'ἄπειρα' in 57d5 (cf. *infra*, note 143). Though he refutes Mugler's interpretation, the issue is not whether it qualifies, as in the classical interpretations, the kinds of each elements, but the existence of a weaker meaning for 'ἄπειρα', not only the 'sense of "infinite," [Mugler] ignoring the possibility that it might have the weaker sense of "indefinitely numerous," or "exceedingly many,"'. In any case, no given number is supposed to limit *a priori*, these different kinds of elements, though there may be only a finite number of them. Vlastos considers probably rightly, it does (cf. *infra*, note 156). For a different (recent) reading strictly limiting the number of kinds of each elements, see O'Brien (1984), p. 342.

[136] Since the faces of the elementary polyhedrons will be multiples of the basic squares or equilateral triangle by unbounded integers. This consequence holds both in the economical and in Timaeus' construction.

[137] Cf. *supra*, a).

[138] The areas of the faces of the different kinds of particles would be all of the form i) for the earth, a square or twice a square a whole number times the basic square i.e. $n^2u^2$ or $2n^2u^2$, ii) a square or thrice a square a whole number times a basic equilateral triangle for the three other elements i.e. $n^2u^2$ or $3n^2u^2$. Thus, the ratios between the different kinds of earth will be squares or twice squares of whole numbers, while the ratios between the different kinds of any of the three other elements (fire, air, water) will be squares or thrice squares of whole numbers. Something certainly worth of an explanation within the mathematical framework of Timaeus' universe.

[139] The sizes of the particles of the different kinds of earth are whole multiples of $u$ and $u\sqrt{2}$ instead of $u$ only, cf. *supra*, c). In modern terms, one considers here the natural situation of the generators on the set of integers $N$. Obviously, if one considers the case on the set $N' = N \cup N\sqrt{2}$, there would be only one generator; however, not only the set $N'$ is much less natural than the set of the integers $N$, but these considerations are coming from the modern set theory, very far from ancient Greek mathematics.

[140] The sizes of the particles of the different kinds of the three other elements (fire, air, water) are whole multiples of $v$ and $v\sqrt{3}$ instead of $v$ only, cf. *supra*, c). Cf. also previous note.

[141] Cf. *supra*, note 80. See also Lloyd (2006), p. 460; Lloyd (2009), p. 13; Artmann-Schäffer (1993), p. 258.

[142] Cf. *supra*, §b).2).

4) The other main difficulty is that there will be many particles of each element, smaller or larger than particles of the other elements. [146] For instance, the basic octahedron has twice the area of the basic tetrahedron. Now, by juxtaposition of triangles, a compound tetrahedron will be at least three times larger than the basic tetrahedron. [147] Thus, necessarily, some particle of fire (tetrahedron) will be larger than the basic particle of air (octahedron). [148] It would then be impossible to maintain the relative size of the particles of the different elements, i.e. that a particle of any kind of fire needs to be smaller than a particle of any kind of air, which in turn needs to be smaller than a particle of any kind of water, [149] although it is how the sensible properties of the different elements are explained by Timaeus. [150]

### ii)    Cornford's solution revisited

It is nevertheless possible to modify Cornford's solution, while avoiding these problems, but keeping the fundamental idea that for each elements, the faces of particles of different kinds are compound of identical basic triangles. [151] Let us first remark that the surfaces, as 2-dimensional figures in a 3-dimensional space, can not only joint but also overlap, to form the faces of the elements, i.e. equilateral triangles and squares. Since Timaeus does not consider

---

[143] 'Consequently, when these are combined amongst themselves and with one another they are infinite in their variety ('διὸ δὴ συμμειγνύμενα αὐτά τε πρὸς αὐτὰ καὶ πρὸς ἄλληλα τὴν ποικιλίαν ἐστὶν ἄπειρα'); and this variety must be kept in view by those who purpose to employ probable reasoning concerning Nature.' (57d4-5).

[144] 56b7-c2; see also for instance, Furley (1987), p. 127.

[145] For a discussion about how large Democritus conceived the size of the atoms, cf. O'Brien (1982).

[146] This problem is already considered in details via the comparison of weights of the different elements in O'Brien (1984), p. 96-99, though not satisfactorily solved (cf. *infra*, notes149 and 150 ).

[147] Cf. *supra*, Figure 7.

[148] Or there would be only one kind of fire, in direct contradiction with Timaeus' account (58c).

[149] 55d8-56b3; according to Denis O'Brien, for the consistency of 'Plato's theory', one has to suppose that 'there is an arbitrary limit on the range of sizes available to any one element, and that the triangular faces of even the largest grade of fire would be smaller, or not much larger, than the triangular faces of even the smallest grade of air or water' (O'Brien (1984), p. 97). Further, Timaeus explains that 'the vessels involved in our sense of smell are too narrow for the varieties of earth and water parts, yet too wide for those of fire and air' (66d4-6) and that 'all odors are finer than water, thicker than air (66e4-5), the different kinds of each elements being considered at once (for some explanation about the ratios involved here, cf. Vlastos (1967), p. 208-209; cf. also Martin (1841), p. 254). The above computations are consequence of the circle of transformations between elements, so it should be stressed that no such restriction applies to the particles (of the different kinds) of earth since this latter element is outside of the cycle of transformations between elements.

[150] For instance, 56a3-4. O'Brien (1984, p. 98-99) agrees (reluctantly) with Martin (1841) on the necessity that the faces of particles of some elements (for instance water) need to be necessarily larger than others (for instance fire). However, under Cornford's construction, with which O'Brien agrees, this would entail for example that 8 basic isosceles triangles will not be able to form a octahedron, but only 2 tetrahedrons; however, 8 times *n* (for some mysterious integer *n*) would be able to form *n* octahedrons (as well as 2*n* tetrahedrons). This claim is not only strange, and without any hint about it in the text, but it is hard to reconcile with Timaeus' account. It is moreover, again, in open contradiction with the laws of transformation (cf. *supra*, beginning of §5). Hence, O'Brien reluctance to this solution is understandable (cf. *ibid.*, note 34, p. 99).

[151] The basic isosceles right triangles and the basic half-equilateral.

directly the faces of the particles but the basic right triangles themselves, which form these faces.[152] The construction for these triangles is done in the two figures below:

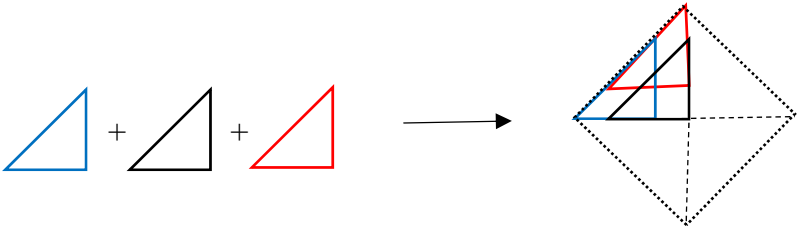

Figure 8 (right isosceles triangles)

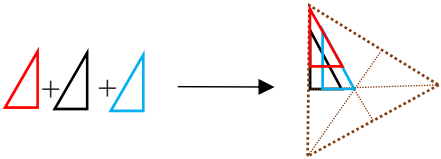

Figure 9 (right scalene triangles)

The first drawing shows how three basic isosceles right triangles form a new larger isosceles right triangle, i.e. one quarter of a larger square. This larger square results from the addition of 12 isosceles basic right triangles. The second one, how three basic half-equilateral right triangles can form a new larger half-equilateral right triangle, six such larger triangles forming a new equilateral triangle. There is an infinite number of such right triangles, isosceles or half-equilateral, all different, and nevertheless all bounded in size. For instance, within the above construction (Figure 8 and Figure 9) using 3 basic right basic triangles, the sides of the largest triangles, either isosceles or half-equilateral, are 3/2 larger than their corresponding basic right triangle,[153] hence, the ratio of their areas is 9/4 (cf. Figure 10 below).[154]

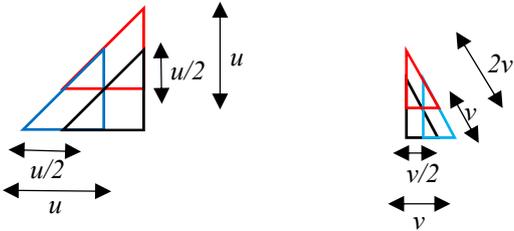

Figure 10

Let us consider among all these right triangles (isosceles or half-equilateral), the ones whose sides are less than 4/3 of the basic right triangles. Then their areas are less than 16/9,[155] hence less than twice the basic right triangles, and there is an 'unlimited' number of such triangles[156]

---

[152] 57c7-d3.
[153] All isosceles right (respectively half-equilateral) triangles are similar, thus the ratios of their corresponding sides are the same.
[154] The ratio of the areas of two similar triangles is the square of the ratio of their sides.
[155] Cf. previous note.
[156] In modern terms, we get all the triangles of sides between the sides of the basic ones and 4/3 of them. Since there is an infinity of different (real) numbers between 1 and 4/3, one obtains an infinity of elementary triangles, forming equilateral triangles of area greater than the basic right triangles (isosceles or half-equilateral) and less than 16/9 their area. Since 16/9 is less than 2, the difficulty of Cornford's construction regarding the relative

all of different sizes. [157]

This construction is not open to the criticism of Cornford's construction with regard to its limited possibilities. [158] Moreover, nor does it have the pitfalls of their own solution, carefully avoided by Cornford. The former authors propose that the elementary right triangles could be arbitrarily small, by halving the basic ones indefinitely. However, there is nothing in Plato's text in favor of such a construction. Moreover, the unity given by the two basic right triangles will be destroyed, and the most important feature of Timaeus' construction as compared to the atomists will be lost. There will be indefinitely different basic triangles, as there are indefinitely different sizes of atoms. Last but not least, it goes against the explicit testimony of Aristotle, [159] so that it is rightly rejected by Cornford. [160] Let us finally remark that, contrary to both Cornford's and Artmann-Schäffer's solutions, the areas of the elementary triangles in our construction are not necessarily commensurable to each other, in agreement with the text, which contains no hint of such relations. For both the isosceles and the scalene basic triangles, it gives an infinite sequence of elementary right triangles of different sizes, commensurable or incommensurable to each other, and nevertheless with areas that are bounded downward and upward. It thus solves two opposing requirements stated by Timaeus: [161]

a) The 'unlimited' complexity of the world is reduced to **two** basic right triangles.

b) The faces of the polyhedrons constructed from these two basic triangles **may** be of **infinite** different sizes, while remaining extremely small and invisible.

Moreover, the relativity of sizes with respect to the four elements, regardless of their kinds, is preserved, in agreement with Timaeus' account. [162] Let us remark that it may seem strange that such a small difference between the basic particles of fire and air (the ratio of their areas is 1 to 2) should result in such different effects. However, as noted both by Plato and by Aristotle, one must consider that even a tiny difference multiplied by a large number may make a huge difference. [163] Such situations are indeed observed in many natural phenomena. [164]

---

sizes of elementary particles (cf. *supra*, §i).b), in particular notes 146 and 150) is avoided. Let us add that this does not mean that all these triangles are to be used, but that the possibilities are 'unlimited' or 'endless', as are the possible different kinds of each element as well (cf. *supra*, note 135). As remarked in note 149, *supra*, the squares are not bound by these restrictions on size, and since they form the faces of particles of (any kind of) earth, neither are the latter, except they must be invisible. Hence, the sides of the isosceles basic right triangle of earth can be much larger than the sides of the basic half-equilateral, so that the particles of earth can be much larger than the particles of the other elements (fire, air, water), thus much more stable.

[157] This is an arbitrary example, for the bounds chosen here, 1 and 4/3, are purely arbitrary. We could have used any magnitude greater than 1 and less than √2 as the upper bound (in order to maintain the relative sizes of the different elementary particles, cf. *supra*, notes 146-156). Moreover, the constraint of size does not apply to the different kinds of water, so that particles of some kinds of water may be much larger than particles of any kind of fire and air, but for reason of stability (cf. *supra*, previous note) smaller than any particle of earth.

[158] See *supra*, §i).b).2); cf. also Pohle (1971), p. 40, Artmann-Schäffer (1993), p. 258.

[159] *On Generation and Corruption*, I, 8, 325b25, cf. *supra*, note 66. Cf. also, *On the Heavens*, III, 7, 30, 306a27-b3.

[160] Timaeus 'stops at a minimum triangle (OBC) of each type, which is taken to be atomic.' (Cornford (1937), p. 234). Nevertheless, some scholars place the property of indefinite divisions of elementary particles under his authority (cf. for example Harte (2002), note 401, p. 240).

[161] 58c5-60e2.

[162] Cf. *supra*, previous paragraph and note 150. This follows from the arbitrary upper bound for all the different kinds of each element (cf. *supra*, note 157).

[163] Cf. *Cratylus*, 436d4-6; Aristotle emphasizes that 'the least initial deviation from the truth is multiplied later a countless of times ('*myrioplasion*') (…) hence that which was small at the start turns out a giant at the end.' (*On the Heavens* I, 5, 271b9-13).

We have tried to show here the amazing simplicity and consistency of Timaeus' universe, founded on mathematical bases, where only two right triangles combine and separate to produce four regular polyhedrons within and outside of a puzzling 'third genus' (the '*khôra*'). The necessary existence of the *khôra*, is a theoretical consequence of the universe as a mathematically ordered set of physical bodies. Its main, if not its only, known property is that it is the 'space' where geometrical planar right triangles (isosceles and half equilateral) get a third dimension as they join together in an orderly way, to form the faces of four regular polyhedrons (the cube, the tetrahedron, the octahedron, and the icosahedron). Conversely, the whole of these polyhedrons is the *khôra*, though not as a static thing, but as a moving sets of bodies. [165] Hence, in a stunning 'tour de force', and in spite of its radical novelty, Timaeus managed to link this geometrical construction to the common physics of his time, founded on the four-element theory, by associating each of these polyhedrons to one of the four elements.

---

[164] For instance, though the difference of mass between two isotopes of uranium, the 298 and 295 is less than 1%, only the second one may be used as the main ingredient to produce atomic bombs.

[165] The dynamical aspect of the *khôra* will be studied in a next article.

# REFERENCES

Translations of the *Timaeus* are quoted from Zeyl (2000) with some slight changes, except for the passage 31b4-32c2, of which a new translation is given *supra*, in the second paragraph.

Translations of Aristotle's works are quoted from the *Works of Aristotle*, edited by W.D. Ross. In particular, his treatise from *On the Heavens*, from J. Stocks' translation, Clarendon Press, 1930 with again some slight changes.